\documentclass[12pt,a4paper]{article}

\usepackage{amsmath,amsfonts}
\usepackage{mathrsfs}

\newtheorem{proposition}{Proposition}

\let\scr\mathscr

\def\zs#1{_{\lower 3pt \hbox{$\scriptstyle#1$}}}
\def\Pb{\mathbf{P}}

\def\Ex{\mathbf{E}}
\def\RR{\mathbb{R}}
\def\sgn{{\rm sgn}} 

\begin{document}
\title{Goodness-of-Fit Tests  for Perturbed Dynamical Systems}
\author{Yury A. \textsc{Kutoyants}\\
{\small Laboratoire de Statistique et Processus, Universit\'e du Maine}\\
{\small 72085 Le Mans, C\'edex 9, France}}
\date{}

\maketitle
%{Running title : }
\begin{abstract}
We consider the goodness of fit testing problem for stochastic differential equation
with small diffusion coefficient. The basic hypothesis is always simple and it
is described by the known trend coefficient. We propose several tests of the
type of Cramer-von Mises, Kolmogorov-Smirnov and Chi-Square. The power
functions of these tests  we study for a special classes of close
alternatives. We discuss the construction of the goodness of fit test based on
the 
local time and the possibility of the construction of asymptotically
distribution free tests in the case of composite basic hypothesis. 
\end{abstract}


\section{Introduction}
We consider the construction of the {\sl goodness-of-fit} (GoF) tests for
dynamical system with small noise, i.e., the observations $X^\varepsilon
=\left\{X_t, 0\leq t\leq T\right\}$ are from the homogeneous stochastic
differential equation
\begin{equation}
\label{1}
{\rm d}X_t=S\left(X_t\right)\;{\rm d}t +\varepsilon \;\sigma
\left(X_t\right)\;{\rm d}W_t,\qquad X_0=x_0, \quad 
0\leq t\leq T 
\end{equation}
with deterministic initial value $x_0$ and known diffusion coefficient $\varepsilon ^2\sigma
\left(\cdot \right)^2>0$. All statistical inference concerns the trend
coefficient $S\left(\cdot \right)$ only. We have two hypotheses: the basic
hypothesis in our consideration is  always  simple
$$
{\scr H}_0:\qquad \quad\qquad \quad S\left(\cdot \right)=S_0\left(\cdot \right)
$$
and the alternative corresponds to the process \eqref{1} with a different
trend coefficient $S\left(\cdot \right)\neq S_0\left(\cdot \right) $. As usual
in GoF testing there are two problems. The first one is to find the threshold
which provides the asymptotic size $\alpha $ of the test and the second is to
describe the behavior of the power function for some classes of
alternatives. There are different ways to present such (nonparametric)
alternatives and we discuss below the choices of alternatives. The problem
considered corresponds to the asymptotics of the small noise, i.e., we study
the properties of the tests as $\varepsilon \rightarrow 0$. We suppose that
the trend and diffusion coefficients satisfy the Lipshits condition
\begin{equation}
\label{2}
\left|S_0\left(x\right)-S_0\left(y\right)\right|+\left|\sigma \left(x\right)-\sigma
\left(y\right)\right| \leq L\;\left|x-y\right|
\end{equation}
hence the equation \eqref{1} has a unique strong solution \cite{LS-01} and we have the
estimates  
\begin{equation}
\label{3}
\Pb_0\left\{\sup_{0\leq t\leq T}\left|X_t-x_t\right|>\delta \right\}\leq
Ce^{-\frac{c}{\varepsilon ^2}},\qquad 
\sup_{0\leq t\leq T}\Ex_0 \left|X_t-x_t\right|^2\leq C\varepsilon ^2,
\end{equation} 
where $C>0,c>0$ are some generic constants (see \cite{FW}, \cite{Kut94}). Here
 the function $\left\{x_t,0\leq t\leq T\right\}$ is solution of the limit
 (deterministic) ordinary equation 
\begin{equation}
\label{4}
\frac{{\rm d}x_t}{{\rm d}t}=S_0\left(x_t\right),\qquad x_0,\quad 0\leq t\leq T.
\end{equation}

Our goal is to present some goodness-of-fit tests for this stochastic model
which are similar to the well-known in classical statistics Cramer-von Mises (C-vM),
Kolmogorov-Smirnov (K-S) and Chi-Square (Ch-S) tests. Remind that these classical tests are
{\sl distribution free}, i.e., their limit distributions do not depend on the
basic hypothesis and therefore the problem of the choice of the threshold is
universal for all tests of such types. Moreover these tests are consistent
against any fixed alternative.  The GoF tests proposed below for the model
\eqref{1} have the similar properties.

 Let us recall the basic properties of the classical tests.  Suppose that we
observe $n$ independent identically distributed random variables
$\left(X_1,\ldots,X_n\right)=X^n$ with continuous distribution function
$F\left(x\right)$ and the basic hypothesis is simple :
$$
{\scr H}_0:\qquad \quad\qquad \quad F\left(x\right)= F_0\left(x\right),\qquad
x\in \RR.
$$
 Then the Cram\'er-von Mises $W_n^2$ and Kolmogorov-Smirnov
$D_n$ statistics are
$$
W_n^2=n\int_{-\infty }^{\infty }\left[\hat
F_n\left(x\right)-F_0\left(x\right)\right]^2 \,{\rm d}F_0\left(x\right),\qquad
D_n=\sup_x\sqrt{n}\left|\hat F_n\left(x\right)-F_0\left(x\right) \right|
$$
respectively. Here
$$
\hat F_n\left(x\right)=\frac{1}{n}\sum_{j=1}^{n}1_{\left\{X_j<x\right\}}
$$
is the empirical distribution function.  Let us denote by
$\left\{W_0\left(s\right), 0\leq s\leq 1\right\}$  a Brownian bridge, i.e., a
continuous Gaussian
process with
$$
\Ex W_0\left(s\right)=0,\qquad \Ex W_0\left(s\right)W_0\left(t\right)= t\wedge
s-st.
$$
Then the
limit behavior of these statistics  can be described with the help of this
process as follows
$$
W_n^2\Longrightarrow \int_{0}^{1}W_0\left(s\right)^2{\rm d}s,\qquad \quad
D_n \Longrightarrow \sup_{0\leq s\leq
1}\left|W_0\left(s\right)\right|.
$$
Hence the corresponding C-vM and K-S tests
$$
\psi_n\left(X^n\right)=1_{\left\{W _n^2>c_\alpha \right\}},\qquad
\phi_n\left(X^n\right)=1_{\left\{D_n>d_\alpha \right\}}
$$
with constants $c_\alpha , d_\alpha  $ defined by the equations
$$
\Pb\left\{\int_{0}^{1}W_0\left(s\right)^2{\rm d}s > c_\alpha \right\}=\alpha
,\qquad \Pb\left\{\sup_{0\leq s\leq 1}\left|W_0\left(s\right)\right|> d_\alpha
\right\}=\alpha
$$
are of asymptotic size $\alpha $.  We see that these tests are
distribution-free (the limit distributions do not depend on the function
$F_*\left(\cdot \right)$) and are consistent against any fixed alternative
(see, for example, Durbin \cite{Dur}, Lehmann and Romano \cite{LR}).

It is interesting to study these tests for {\sl non degenerate set of
alternatives}, i.e., for alternatives with limit power function is grater than
$\alpha $ and less than 
1. It can be realized on the close nonparametric alternatives of the special
form making this problem asymptotically equivalent to the {\sl signal in
Gaussian noise} problem. Let us put
$$
F\left(x\right)=F_0\left(x\right)+\frac{1}{\sqrt{n}}\int_{-\infty
}^{x}h\left(F_0\left(y\right)\right)\;{\rm d}F_0\left(y\right) ,
$$
where the function $h\left(\cdot \right)$ describes the alternatives.
We suppose that
$$
\int_{0}^{1}h\left(s\right)\;{\rm d}s=0,\qquad
\int_{0}^{1}h\left(s\right)^2\;{\rm d}s<\infty .
$$
Then we
have the following convergence (under fixed alternative, given by the function
$h\left(\cdot \right)$):
\begin{align*}
&W_n^2\Longrightarrow \int_{0}^{1}\left[\int_{0}^{s}h\left(v\right){\rm
d}v+W_0\left(s\right)\right]^2{\rm d}s ,\\
& D_n \Longrightarrow
\sup_{0\leq s\leq 1}\left|\int_{0}^{s}h\left(v\right){\rm
d}v+W_0\left(s\right)\right|
\end{align*}
We see that this problem is asymptotically equivalent to the  following {\sl
signal in Gaussian noise} problem:
\begin{equation}
\label{5}
{\rm d}Y_s=h_*\left(s\right)\,{\rm d}s+{\rm d}W_0\left(s\right),\quad 0\leq
s\leq 1.
\end{equation}

Indeed, if we use the statistics
$$
W^2=\int_{0}^{1}Y_s^2\;{\rm d}s,\qquad D=\sup_{0\leq s\leq 1}\left|Y_s\right|
$$
then under hypothesis $h\left(\cdot \right)\equiv 0$ and alternative
$h\left(\cdot \right)\neq  0$   the distributions of these
statistics coincide with the limit distributions of $W_n^2$ and
$D_n$ under hypothesis and alternative respectively.

\section{C-vM and K-S type tests}

\subsection{Choice of the thresholds}
We test the basic simple  hypothesis ${\scr H}_0$ and our goal is to
study the GoF tests of asymptotic size $\alpha \in\left(0,1\right)$. Let us
denote the class of such tests as 
$$
{\cal K}_\alpha =\left\{\psi_\varepsilon\; : \quad \Ex_0 \psi_\varepsilon\left(X^\varepsilon \right)=\alpha +o\left(1\right)\right\},
$$
where $\Ex_0$ is expectation under hypothesis. 
We use the following regularity condition.

\bigskip

{\bf Condition ${\cal R}$.} {\sl The function $S_0\left(\cdot \right)$ has
two continuous bounded derivatives $ S_0'\left(x\right),
S_0''\left(x\right) $ and the function  
$\sigma \left(x\right)$ has one 
continuous bounded derivative $  \sigma' \left(x\right)$ and the both
functions are positive :
$S_0\left(x\right)>0, \sigma^2 \left(x\right)>0$  for  $x\geq x_0$ }. 

\bigskip

Remind that under this condition the equation \eqref{1} has a unique strong
solution and this solution converges uniformly on $t\in \left[0,T\right]$ to
the deterministic solution of the equation \eqref{4} (see, e.g., \cite{LS-01}
\cite{Kut94}). Moreover, the stochastic process $X_t$ is differentiable
w.r.t. $\varepsilon $ at the point $\varepsilon =0$ and its derivative
$x_t^{\left(1\right)}$ satisfies the linear equation
\begin{equation}
\label{5a}
{\rm d}x_t^{\left(1\right)}= S_0'\left(x_t\right)x_t^{\left(1\right)}\,{\rm
d}t+ \sigma \left(x_t\right)\,{\rm d}W_t,\quad x_0^{\left(1\right)}=0. 
\end{equation}
 For the proof see, e.g., \cite{Kut94}, Lemma 3.3.

 To construct the C-vM and K-S type tests we use
the statistics
$$
\delta _\varepsilon =\left[\int_{0}^{T}\left(\frac{\sigma
\left(x_t\right)}{S_0\left(x_t\right)}\right)^{2}{\rm 
d}t\right]^{-2} \; \int_{0}^{T}\left(\frac{X_t-x_t}{\varepsilon
\,S_0\left(x_t\right)^{2}}\right)^2\;\sigma \left(x_t\right)^2\,{\rm d}t. 
$$
and
$$
\gamma _\varepsilon =\left[\int_{0}^{T}\left(\frac{\sigma
\left(x_t\right)}{S_0\left(x_t\right)}\right)^{2}{\rm 
d}t\right]^{-1/2} \; \sup_{0\leq t\leq T}\left|\frac{X_t-x_t}{\varepsilon
\,S_0\left(x_t\right)}\right| 
$$
respectively.

Below $c_\alpha $ and $b_\alpha $ are solutions of the equations

\begin{equation}
\label{6}
\Pb\left\{\int_{0}^{1}w_v^2\;{\rm d}v>c_\alpha \right\}=\alpha,\qquad \quad 
\Pb\left\{\sup_{0\leq v\leq 1}\left|w_v\right|>b_\alpha \right\}=\alpha,
\end{equation}
where $w_v,0\leq v\leq 1$ is some Wiener process. 

\begin{proposition}
\label{P1} Let the condition ${\cal R}$ be fulfilled then the tests
$\psi_\varepsilon  =1\zs{\left\{\delta _\varepsilon >c_\alpha \right\}}$ and
$\phi_\varepsilon  =1\zs{\left\{\gamma  _\varepsilon >b_\alpha \right\}}$
belong to the class ${\cal K}_\alpha $.
\end{proposition}
{\bf Proof.} The stochastic process $\varepsilon ^{-1}\left(X_t-x_t\right)$
converges in probability uniformly on $t\in \left[0,T\right]$ to the limit
$x_t^{\left(1\right)}$ - solution of the linear equation \eqref{5a}, i.e., for
any $\kappa >0$
$$
\Pb_0\left\{\sup_{0\leq t\leq T}\left|\frac{X_t-x_t}{\varepsilon
}-x_t^{\left(1\right)}\right|>\kappa \right\} \longrightarrow 0. 
$$
This solution can be written explicitly as
$$
x_t^{\left(1\right)}=\int_{0}^{t}\exp\left\{\int_{s}^{t}
S_0'\left(x_v\right){\rm d}v\right\}\sigma \left(x_s\right)\;{\rm d}W_s. 
$$
Below we follow \cite{Kut94}, where the similar calculus were done.  Using
\eqref{2} we can write
\begin{align}
\int_{s}^{t} S_0'\left(x_v\right){\rm d}v&=\int_{s}^{t}\frac{
S_0'\left(x_v\right) }{S_0\left(x_v\right)}\;\frac{{\rm d}x_v}{{\rm d}v}\; {\rm
d}v=\int_{s}^{t}\frac{
S_0'\left(x_v\right) }{S_0\left(x_v\right)}\; {\rm d}x_v=\nonumber\\
&=\int_{x_s}^{x_t}\frac{
S_0'\left(x\right) }{S_0\left(x\right)}\; {\rm d}x=\int_{x_s}^{x_t} 
\Bigl(\ln  S_0\left(x\right)\Bigr)' \; {\rm d}x=\ln
\frac{S_0\left(x_t\right)}{S_0\left(x_s\right) }.
\label{cal}
\end{align}
Hence 
$$
x_t^{\left(1\right)}=S_0\left(x_t\right)\;\int_{0}^{t}\frac{\sigma
\left(x_s\right) }{S_0\left(x_s\right)}\;{\rm 
d}W_s=S_0\left(x_t\right)\;W\left( \int_{0}^{t}
\frac{\sigma \left(x_s\right)^2 }{S_0\left(x_s\right)^2} \;{\rm d}s\right) 
$$
where $W\left(\cdot \right)$ is some Wiener process. Further
\begin{align*}
\int_{0}^{T}\left(\frac{X_t-x_t}{\varepsilon
\,S_0\left(x_t\right)^{2}}\right)^2\sigma \left(x_s\right)^2{\rm d}t &\rightarrow
\int_{0}^{T}W\left( \int_{0}^{t} 
\frac{\sigma \left(x_s\right)^2}{S_0\left(x_s\right)^2}{\rm d}s  \right) {\rm
d}\left(\int_{0}^{t}\frac{\sigma
\left(x_s\right)^2}{S_0\left(x_s\right)^2}{\rm d}s\right)\\ 
& =\int_{0}^{u _T}W\left(u\right)^2\;{\rm d}u =u _T^2\int_{0}^{1}w_v^2\;{\rm d}v ,
\end{align*}
where we put 
$$
u=\int_{0}^{t}\frac{\sigma \left(x_s\right)^2}{S_0\left(x_s\right)^2}\;{\rm
d}s=v\;\int_{0}^{T}\frac{\sigma
\left(x_s\right)^2}{S_0\left(x_s\right)^2}\;{\rm 
d}s =v\;u _T.
$$
Here $w_v=u_T^{-1/2}W\left(u_Tv\right),0\leq v\leq 1$ is a Wiener process.

Hence (under hypothesis ${\scr H}_0$) 
$$
\delta _\varepsilon \longrightarrow \int_{0}^{1}w_v^2\;{\rm d}v .
$$
For the statistic $\gamma _\varepsilon $ the similar consideration leads to
the relation
$$
\gamma _\varepsilon \longrightarrow \sup_{0\leq v\leq 1}\left|w_v\right|.
$$
Therefore  for the first type errors we have
$$
\alpha _\varepsilon\left(\delta _\varepsilon \right) =\Pb_0\left\{\delta _\varepsilon >c_\alpha
\right\}\longrightarrow \Pb\left\{\int_{0}^{1}w_v^2\;{\rm d}v>c_\alpha  \right\}=\alpha, 
$$
and
$$
\alpha _\varepsilon\left(\gamma _\varepsilon \right) =\Pb_0\left\{\gamma  _\varepsilon >b_\alpha
\right\}\longrightarrow \Pb\left\{ \sup_{0\leq v\leq
1}\left|w_v\right|>b_\alpha \right\}=\alpha.  
$$
The Proposition \ref{P1} is proved. 

\bigskip

\subsection{Similar tests.}
Note that the similar result we have if we use  the statistics
$$
\bar\delta _\varepsilon =\left[\int_{0}^{T}\left(\frac{\sigma
\left(X_t\right)}{S_0\left(X_t\right)}\right)^{2}{\rm 
d}t\right]^{-2} \; \int_{0}^{T}\left(\frac{X_t-x_t}{\varepsilon
\,S_0\left(X_t\right)^{2}}\right)^2\;\sigma \left(X_t\right)^2\,{\rm d}t 
$$
and
$$
\bar\gamma _\varepsilon =\left[\int_{0}^{T}\left(\frac{\sigma
\left(X_t\right)}{S_0\left(X_t\right)}\right)^{2}{\rm 
d}t\right]^{-1/2} \; \sup_{0\leq t\leq T}\left|\frac{X_t-x_t}{\varepsilon
\,S_0\left(X_t\right)}\right|. 
$$
because, as we mentioned above, the process $X_t$ converges uniformly on $t\in
\left[0,T\right]$ to the deterministic solution $x_t$ and this implies the
convergence
$$
\bar\delta _\varepsilon \longrightarrow \int_{0}^{1}w_v^2\,{\rm
d}v,\qquad \qquad \bar\gamma_\varepsilon \longrightarrow \sup_{0\leq v\leq 1}\left|w_v\right|.
$$

\bigskip

These tests can be slightly  simplified, if we replace the equation \eqref{4} by the
following one
\begin{equation}
\label{}
\frac{{\rm d}\hat X_t}{{\rm d}t}=S_0\left(X_t\right),\qquad \hat X_0=x_0,
\end{equation}
i.e., we put
$$
\hat X_t=x_0+\int_{0}^{t}S_0\left(X_s\right)\;{\rm d}s,\qquad \quad 0\leq t\leq T,
$$
(we need not to solve the equation \eqref{4}, just to calculate the integral,
using observations) and introduce the statistic
$$
\hat\delta  _\varepsilon =\left[\int_{0}^{T}\sigma \left(X_t\right)^2{\rm
d}t\right]^{-2} \;\int_{0}^{T}\sigma \left(X_t\right)^2\left(\frac{X_t-\hat
X_t}{\varepsilon }\right)^2 \;{\rm d}t 
$$
Then under hypothesis ${\scr H}_0$
\begin{align*}
&\varepsilon ^{-2}\int_{0}^{T}\sigma \left(X_t\right)^2\left(X_t-\hat X_t\right)^2 \;{\rm
d}t=\int_{0}^{T}\sigma \left(X_t\right)^2\left[\int_{0}^{t}\sigma 
\left(X_s\right){\rm d}W_s\right]^2{\rm d}t \\
&\qquad =\int_{0}^{T}\left[W\left(\int_{0}^{t}\sigma 
\left(X_s\right)^2{\rm d}s\right)\right]^2{\rm d}\left(\int_{0}^{t}\sigma
\left(X_s\right)^2{\rm d}s\right)  \\ 
&\qquad =\int_{0}^{\tau_T }W\left(r\right)^2\,{\rm
d}r\longrightarrow\int_{0}^{\tau_T^o }W\left(r\right)^2\,{\rm
d}r\ =\left(\tau_T^o\right)^2\int_{0}^{1}w_v^2\,{\rm d}v, 
\end{align*}
where 
$$
\tau _T=\int_{0}^{T}\sigma \left(X_s\right)^2{\rm d}s\longrightarrow \tau_T^o=\int_{0}^{T}\sigma \left(x_s\right)^2{\rm d}s,\qquad \quad
    v=\frac{r}{\tau _T^o}.
$$
Hence we have the convergence 
$$
\hat\delta  _\varepsilon \Longrightarrow \int_{0}^{1}w_v^2\,{\rm d}v.
$$
Of course, we have the {\sl distribution free} limit for the corresponding
statistic $\hat \gamma _\varepsilon $ too.

\subsection{Partially observed linear system}

Suppose that we observe a random process $X^T=\left\{X_t,0\leq
t\leq T\right\}$ and we have to test the hypothesis ${\scr H}_0$ that this
process comes from  the following linear partially observed  system 
\begin{align*}
{\rm d}Y_t&=A_t\,Y_t\,{\rm d}t+\varepsilon\, B_t\,{\rm d}V_t,\quad
Y_0=y_0,\quad 0\leq t\leq T,\\
{\rm d}X_t&=C_t\,Y_t\,{\rm d}t+\varepsilon\sigma _t\, {\rm d}W_t,\quad X_0=0,\quad
0\leq t\leq T,
\end{align*}
where $A_t,B_t$,  $C_t$ and  $\sigma _t$  are known functions and $V_t$ and
$W_t$ are 
independent Wiener processes. We suppose as well that these functions satisfy
the {\it usual conditions} which allows us to write the equations of
filtration (see Liptser and Shiryaev \cite{LS-01}). 

If the hypothesis is true then according to well known Kalman-Bucy theory (see
Liptser and Shiryaev \cite{LS-01}, Theorem 10.1)  the conditional expectation
$M_t=\Ex_0 
\left(Y_t|X_s,0\leq s\leq t \right)$ satisfies the equations
\begin{align*}
{\rm d}M_t&=A_t\,M_t\,{\rm d}t+\frac{C_t\,\gamma \left(t\right)}{\varepsilon
^2\sigma _t^2} \left[{\rm d}X_t-C_t\,M_t\,{\rm d}t\right],\\
\frac{{\rm d}\gamma \left(t\right)}{{\rm d}t}&=2A_t\gamma
\left(t\right)-\frac{C_t^2\,\gamma \left(t\right)^2}{\varepsilon 
^2\sigma _t^2}+\varepsilon ^2\,B_t^2
\end{align*}
with initial values $M_0=\Ex_0 Y_0=y_0$ and $\gamma
\left(0\right)=\Ex_0\left(Y_0-\Ex_0 Y_0\right)^2=0$. Note that if we put
$\Gamma_t=\gamma \left(t\right)\varepsilon ^{-2} $, then this system can be
rewritten as 
 \begin{align*}
{\rm d}M_t&=A_t\,M_t\,{\rm d}t+\frac{C_t\,\Gamma_t}{\sigma _t^2} \left[{\rm
d}X_t-C_t\,M_t\,{\rm d}t\right],\quad M_0=y_0,\\ 
\frac{{\rm d}\Gamma_t \left(t\right)}{{\rm d}t}&=2A_t\Gamma_t
-\frac{C_t^2\,\Gamma_t^2}{\sigma _t^2}+B_t^2,\qquad \Gamma_0=0.
\end{align*}
Remind as well that the observed process admits the representation
$$
{\rm d}X_t=C_t\,M_t\,{\rm d}t+\varepsilon \sigma _t\,{\rm d}\bar W_t,\quad
0\leq t\leq T, 
$$
where $\bar W_t$ is {\it innovation} Wiener process defined by this equality. 
This representation suggests to define the statistic
$$
\delta _\varepsilon =\left(\varepsilon  \int_{0}^{T}\sigma
\left(t\right)^2{\rm d}t \right)^{-2}\int_{0}^{T}\sigma
\left(t\right)^2\left[X_t-\int_{0}^{t}C_s\,M_s\,{\rm d}s\right]^2 {\rm d}t.
$$
Elementary calculations yield
\begin{align*}
&\varepsilon ^{-2}\int_{0}^{T}\sigma
\left(t\right)^2\left[X_t-\int_{0}^{t}C_s\,M_s\,{\rm d}s\right]^2 {\rm d}t
=\int_{0}^{T}\sigma 
\left(t\right)^2\left[\int_{0}^{t}\sigma _s\,{\rm d}\bar W_s\right]^2 {\rm
d}t\\
&\qquad =\int_{0}^{T}\sigma 
\left(t\right)^2\left[\bar W\left(\int_{0}^{t}\sigma _s^2\,{\rm
d}s\right)\right]^2 {\rm d}t=\left(\int_{0}^{T}\sigma 
\left(t\right)^2\, {\rm d}t\right)^2\int_{0}^{1}w_v^2{\rm d}v .
\end{align*}
Hence the statistic $\delta _\varepsilon $ (under hypothesis) is
$$
\delta _\varepsilon =\int_{0}^{1}w_v^2\;{\rm d}v
$$
and the corresponding test $\psi_\varepsilon =1_{\left\{\delta _\varepsilon
>c_\alpha \right\}}$ is distribution free.

\subsection{Alternatives}

Let us consider nonparametric alternatives $S\left(\cdot \right)\neq
S_0\left(\cdot \right)$, which correspond to the equation
\begin{equation*}
{\rm d}X_t=S\left(X_t\right)\,{\rm d}t+\varepsilon \sigma
\left(X_t\right)\;{\rm d}W_t,\quad X_0=x_0,\quad 0\leq t\leq T.
\end{equation*}
 We suppose that the functions $S\left(\cdot \right)$ and $\sigma \left(\cdot
\right)$ satisfy the 
conditions,
\begin{align*}
&\left|S\left(x\right)-S\left(y\right)\right|+\left|\sigma
\left(x\right)-\sigma\left(y\right)\right|\leq L\left|x-y\right|,\\
&\left|S\left(x\right)\right|+\left|\sigma
\left(x\right)\right|\leq L\left(1+\left|x\right|\right).
\end{align*}
 The set of such functions we denote as ${\cal S}\left(L\right)$ and the limit
of the stochastic process $X_t,0\leq t\leq T$ we write as $x_t\left(S\right),
0\leq t\leq T$.

There are many ways to introduce the class of alternatives. We consider two of
them. Let us define the sets
\begin{align*}
{\cal F}_r&=\left\{S\left(\cdot \right)\in{\cal S}\left(L\right):\;\;
\left\|x_.\left(S\right)-x_.\right\|\geq r\right\},\\ 
{\cal G}_r&=\left\{S\left(\cdot \right)\in{\cal S}\left(L\right):\;\;
\left\|S\left(x_.\right)-S_0\left(x_.\right)\right\|\geq r\right\}, 
\end{align*}
where $\left\|\cdot \right\|$ is ${\cal L}_2\left(0,T\right)$-norm, say,
$$
\left\|S\left(x_.\right)-S_0\left(x_.\right)\right\|^2=\int_{0}^{T}
\left[S\left(x_t\right)-S_0\left(x_t\right)\right]^2{\rm d}t.
$$
Here $x_t,\leq 0\leq t\leq T$ is solution of the equation \eqref{4} (under
hypothesis ${\scr H}_0$).

Let us start with the problem
\begin{align*}
{\scr H}_0\;\;&:\qquad \quad S\left(\cdot \right)=S_0\left(\cdot \right),\\
{\scr H}_1\;\;&:\qquad \quad S\left(\cdot \right)\in {\cal F}_r
\end{align*}
To show the consistency of the test $\psi_\varepsilon $ note that $\delta
_\varepsilon  \geq \kappa \varepsilon ^{-1}\left\|X_.-x_.\right\|$ with some
$\kappa >0$. We have
\begin{align*}
&\inf_{S\in {\cal F}_r}\Pb_S\left\{\delta _\varepsilon >c_\alpha \right\}\geq
\inf_{S\in {\cal F}_r}\Pb_S\left\{\kappa \varepsilon
^{-1}\left\|X_.-x_.\right\| >\sqrt{c_\alpha} \right\} \\
&\quad \geq  \inf_{S\in {\cal F}_r}\Pb_S\left\{\varepsilon
^{-1}\left\|x_.\left(S\right)-x_.\right\|-\varepsilon
^{-1}\left\|X_.-x_.\left(S\right)\right\| >\sqrt{c_\alpha} \right\} \\
&\quad \geq  1-\sup_{S\in {\cal F}_r}\Pb_S\left\{\varepsilon
^{-1}\left\|X_.-x_.\left(S\right)\right\| >\varepsilon
^{-1}\left\|x_.\left(S\right)-x_.\right\|-\kappa ^{-1}\sqrt{c_\alpha}
\right\}  \\
&\quad \geq  1-\left(\varepsilon
^{-1}r-\kappa ^{-1}\sqrt{c_\alpha}\right)^{-2}\sup_{S\in {\cal
F}_r}\Ex_S\int_{0}^{T}\left( \frac{X_t-x_t\left(S\right)}{\varepsilon
}\right)^2{\rm d}t \geq 1-C\varepsilon ^2.
\end{align*}

Therefore even if $r=r_\varepsilon \rightarrow 0$ such that $\varepsilon
^{-1}r_\varepsilon \rightarrow \infty $, then the C-vM test is {\sl
minimax consistent}, i.e., its power function tends to 1 uniformly on
$S\left(\cdot \right)\in {\cal F}_r $. The same is true for the
Kolmogorov-Smirnov test. It is not so in the case of the hypotheses testing
problem 
\begin{align*}
{\scr H}_0\;\;&:\qquad \quad S\left(\cdot \right)=S_0\left(\cdot \right),\\
{\scr H}_1\;\;&:\qquad \quad S\left(\cdot \right)\in {\cal G}_r
\end{align*}

It will be more convenient to write
$S\left(x\right)-S_0\left(x\right)$ as $\varepsilon h\left(x\right)\sigma 
\left(x\right)^2S_0\left(x\right)^{-1}$, where $h\left(\cdot \right)$ is
such that $S\left(\cdot \right)\in {\cal G}_r $. This corresponds to the model
\begin{equation}
\label{7}
{\rm d}X_t=S_0\left(X_t\right)\,{\rm d}t+\varepsilon \;
\frac{h\left(X_t\right)\sigma \left(X_t\right)^2}{S_0\left(X_t\right)}\;{\rm
d}t+\varepsilon \sigma \left(X_t\right)\;{\rm d}W_t,\quad 0\leq t\leq T
\end{equation}
with the same initial value $X_0=x_0$.   The case
$h\left(\cdot \right)\equiv 0$ corresponds to the hypothesis ${\scr H}_0$. 
We start with the  study of the limit behavior of the statistic $\delta
_\varepsilon $ 
under a fixed  alternative $h\left(\cdot \right)$. 

Let us denote by $x_t^h$ the solution of the equation
$$
\frac{{\rm d}x_t^h}{{\rm d}t }=S_0\left(x_t^h\right) +\varepsilon \;
\frac{\sigma \left(x_t^h\right)^2}{S_0\left(x_t^h\right)}\;h\left(
x_t^h\right),\qquad x_0^h=x_0; 
$$
and write 
$$
\frac{X_t-x_t}{\varepsilon }=\frac{X_t-x_t^h}{\varepsilon
}+\frac{x_t^h-x_t}{\varepsilon } .
$$
It is easy to see that 
$$
\frac{X_t-x_t^h}{\varepsilon}\longrightarrow x_t^{\left(1\right)},\qquad 0\leq
t\leq T 
$$
in probability and the direct calculations yield
\begin{align*}
\frac{x_t^h-x_t}{\varepsilon }&=\int_{0}^{t}\frac{S_0\left(x_s^h
\right)-S_0\left(x_s \right)}{\varepsilon }\;{\rm d}s+
\int_{0}^{t}\frac{\sigma \left(x_s^h\right)^2h\left(x_s^h 
\right)}{S_0\left(x_s^h \right)
}\;{\rm d}s\\
&=\int_{0}^{t}\left(\frac{x_s^h-
x_s }{\varepsilon }\right)\;  S_0'\left(\tilde x_s^h \right){\rm d}s+
\int_{0}^{t}\frac{\sigma \left(x_s^h \right)^2h\left(x_s^h
\right)}{S_0\left(x_s^h \right) }\;{\rm d}s.
\end{align*}
Hence $\frac{x_t^h-x_t}{\varepsilon }$ converges to the function $\dot x_t^h$
which is solution of the equation
$$
\frac{{\rm d}\dot x_t^h}{{\rm d}t}= S_0'\left( x_t \right)\,\dot
x_t^h+\frac{\sigma \left(x_t  \right)^2}{S_0\left(x_t \right) }\;h\left(x_t
\right),\qquad \dot x_0^h=0. 
$$
Therefore
\begin{align*}
\dot x_t^h&=\int_{0}^{t}\exp\left\{\int_{s}^{t} S_0'\left(x_v\right)\;{\rm
d}v\right\} \;\frac{\sigma \left(x_s \right)^2}{S_0\left(x_s \right)
}\;h\left(x_s  \right)\;{\rm d}s\\
& =S_0\left(x_t\right)\;\int_{0}^{t}\frac{\sigma \left(x_s \right)^2}{S_0\left(x_s
\right)^2 }\;h\left(x_s  \right)\;{\rm d}s 
\end{align*}
and
\begin{align*}
\frac{X_t-x_t}{\varepsilon S_0\left(x_t\right)}&\longrightarrow
Z\left(x_t\right)\equiv 
\frac{x_t^{\left(1\right)}+\dot x_t^h}{S_0\left(x_t\right) }\\
&=\int_{0}^{t}\frac{\sigma \left(x_s
\right)^2}{S_0\left(x_s 
\right)^2 }\;h\left(x_s  \right)\;{\rm d}s  +W\left(\int_{0}^{t}\frac{\sigma
\left(x_s\right)^2  }{S_0\left(x_s 
\right)^2 }\;{\rm d}s  \right)   \\
&=\int_{x_0}^{x_t}\frac{\sigma \left(x
\right)^2}{S_0\left(x 
\right)^3 }\;h\left(x  \right)\;{\rm d}x +W\left(\int_{x_0}^{x_t}\frac{\sigma
\left(x\right)^2  }{S_0\left(x 
\right)^3 }\;{\rm d}x \right) .
\end{align*}
Let us put 
$$
u\left(x\right)=\int_{x_0}^{x}\frac{\sigma \left(y
\right)^2}{S_0\left(y\right)^3 }\;{\rm d}y\quad \in \quad \left[0,u_T\right]
$$
and denote $x\left(u\right)$ the function inverse to $u\left(x\right) $
These relations yield the following representation for the limit of the test
statistic
\begin{align}
\delta _\varepsilon \longrightarrow \; \delta _0 & =
\int_{0}^{T}
\;\frac{Z\left(x_t\right) ^2 \sigma \left(x_t\right)^2}{u_T^{2} \;
S_0\left(x_t\right)^{2} }\,{\rm d} t =
\int_{0}^{T}\frac{Z\left(x_t\right)^{2}}{u_T^{2}}  
\;\,{\rm d} \left(\int_{x_0}^{x_t} \frac{\sigma \left(y\right)^2}{
S_0\left(y\right)^{3} }\;{\rm d}y\right) \nonumber\\
&=\int_{x_0}^{x_T}\frac{Z\left(x\right)^{2}}{u_T^{2}}  
\;\,{\rm d} \left(\int_{x_0}^{x} \frac{\sigma \left(y\right)^2}{
S_0\left(y\right)^{3} }\;{\rm d}y\right)\nonumber\\
&=u_T^{-2}\int_{0}^{u_T}\left[\int_{0}^{u}h\left(x\left(z\right)\right){\rm
d}z+W\left(u\right) \right]^2 {\rm d}u \nonumber\\
&=\int_{0}^{1}\left[\int_{0}^{v}h_*\left(s\right){\rm
d}s+w_v \right]^2 {\rm d}s,
\label{8}
\end{align}
where we put $u=u_T\,v$ and denoted
$h_*\left(s\right)=u_T^{1/2}h\left(x\left(u_Ts\right)\right)$ and
$w_v=u_T^{-1/2} W\left(u_Tv\right)$. It is easy to see that $w_v,0\leq v\leq
1$ is standard Wiener process. 

Remind that if the observed process is of type {\sl signal in white Gaussian
noise}:
$$
{\rm d}Y_s=h\left(s\right)\,{\rm d}s+{\rm d}w_s,\quad 0\leq
s\leq 1, 
$$ 
then the natural distance between hypothesis $h\left(s\right)\equiv 0$ and
$h\left(s\right)\neq 0$ is ${\cal L}^2\left(0,1\right)$:
$$
\int_{0}^{1}h\left(s\right)^2\;{\rm d}s\geq \rho^2 .
$$ 

In our case it corresponds to 
$$
\int_{0}^{1}h_*\left(s\right)^2\;{\rm
d}s=\int_{x_0}^{x_T}h\left(x\right)^2\frac{\sigma
\left(x\right)^2}{S_0\left(x\right)^3}\;{\rm d}x \geq \rho^2.
$$

The convergence \eqref{8} provides  the limit of the power function given in
the next proposition. 
\begin{proposition}
\label{P2} let the condition ${\cal R}$ be fulfilled, then for any  function
$h\left(\cdot \right)\in {\cal H}_r$ we have
\begin{align*}
\beta \left(\delta _\varepsilon ,h\right)=\Pb_h\left\{\delta _\varepsilon
>c_\alpha \right\}=
\Pb\left\{\int_{0}^{1}\left[\int_{0}^{v}h_*\left(s\right){\rm 
d}s+w_v \right]^2 {\rm d}s>c_\alpha  \right\}+o\left(1\right).  
\end{align*}
\end{proposition}

\bigskip

Let us define the
composite (nonparametric) alternative  as
\begin{align*}
&{\scr H}_1:\quad \quad   &h\left(\cdot \right)\in 
{\cal H}_{r },
\end{align*} 
where
$$
{\cal H}_{r}=\left\{h\left(\cdot \right)\; : \qquad \left\|\frac{h\left(x_.
\right)\sigma \left(x_.\right)^2}{S_0\left(x_.\right)}\right\| \geq
\frac{r}{\varepsilon } 
 \right\} .
%\quad \qquad \left\|h\left(\cdot
%\right)\right\|_\mu^2=\int_{x_0}^{x_T}h\left(x\right)^2 \;\mu \left({\rm
%d}x\right). 
$$
%The ``natural choice'' of the measure $\mu \left(\cdot \right)$ will be done later. 
Note that $x_.$ here is ``from hypothesis ${\scr H}_{0}$''. As $h\left(\cdot
\right)$ is an arbitrary function, this alternative coincides with ${\cal
G}_r$.  Let us show that without regularity conditions the problem of
hypotheses testing is degenerate in the following sens: we have for the power
of the test
$$
\inf_{S\in {\cal G}_{r } }\beta \left(\delta _\varepsilon
,h\right) \longrightarrow \alpha 
$$
as $\varepsilon \rightarrow 0$ even for fixed $r>0$. The condition
$S\left(\cdot \right)\in {\cal S}\left(L\right)$ already provides the control
of the first derivative, that is why we have to  weaken it slightly. At
particularly,  we suppose that
$S\left(\cdot \right)\in {\cal S}\left(L\right)$, but $L>0$ can be  not  the
same for different functions $S\left(\cdot \right)$. 

Let us put 
$$
h_n\left(x\right)=c\;\frac{S_0\left(x\right)^2}{\sigma \left(x\right)^2}\;\cos
\left(n\left(x-x_0\right)\right) 
$$
then,
$$
\left|S_n\left(x\right)-S_0\left(x\right)\right|=\varepsilon
\,c\,\left|S_0\left(x\right)\;\cos
\left(n\left(x-x_0\right)\right) \right|
$$
and $L=K\,c\,\varepsilon \,n$ with some $K>0$. Further, using $2\cos^2\varphi
 =1-\cos \left(2\varphi \right)$ we obtain
\begin{align*}
&\int_{x_0}^{x_T}h_n\left(x\right)^2\frac{\sigma
\left(x\right)^2}{S_0\left(x\right)^3}\;{\rm d}x
=\frac{c^2}{2}\int_{x_0}^{x_T}\frac{S_0\left(x\right)}{\sigma
\left(x\right)^2}\; {\rm
d}x\\
&\qquad  -\frac{c^2}{2}\int_{x_0}^{x_T}\cos\left(2n\left(x-x_0\right)\right)\; 
\frac{S_0\left(x\right)}{\sigma \left(x\right)^2}\; {\rm d}x\longrightarrow
\frac{c^2}{2}\int_{x_0}^{x_T}\frac{S_0\left(x\right)}{\sigma 
\left(x\right)^2}\; {\rm d}x
\end{align*}
as $n\rightarrow \infty $.  Remind that $S_0\left(x\right)$ and $\sigma
\left(x\right)$ are continuous positive  functions. 
The constant $c=c\left(\varepsilon \right)>0$ can be chosen from the condition
$$
\frac{c^2}{4}\int_{x_0}^{x_T}\frac{S_0\left(x\right)}{\sigma 
\left(x\right)^2}\; {\rm d}x\geq \frac{r^2}{\varepsilon ^2} .
$$
Hence $h_n\left(\cdot \right)\in {\cal H}_r $. On the other hand
\begin{align*}
\int_{x_0}^{x}h_n\left(y\right)\,\frac{\sigma \left(y\right)^2}{S_0
\left(y\right)^3}\;{\rm d}y &=\frac{c\left(\varepsilon
\right)}{n}\frac{\sin\left(n\left(x-x_0\right)\right)}{ S_0 \left(x\right) }\\
&+\frac{c\left(\varepsilon \right)}{n} \int_{x_0}^{x}\frac{
S_0'\left(y\right)\sin\left(n\left(y-x_0\right)\right)}{S_0\left(y\right)^2
}{\rm d}y \longrightarrow 0
\end{align*}
uniformly in $x\in \left[x_0,x_T\right]$ if we put, say, $n=n\left(\varepsilon
\right)=c\left(\varepsilon \right)^2$. Therefore 
$$
\lim_{\varepsilon \rightarrow 0}\inf_{h\left(\cdot \right)\in {\cal H}_r
}\beta \left(\delta _\varepsilon 
,h\right)\leq \lim_{\varepsilon \rightarrow 0} \inf_{h_n\left(\cdot \right)
}\beta \left(\delta _\varepsilon 
,h_n\right)=\alpha. 
$$
Hence if we use the introduced above statistics then this hypotheses testing
problem will be asymptotically degenerated. It can be shown that this is not
the particular property of these two tests, but is true for any other tests
too.

\section{Chi Square Test}

If the alternative is defined by the inequality
$$
\left\|x_.\left(S\right)-x_.\right\|\geq r,
$$
then it is natural to replace $x_t\left(S\right)$ by its estimate $X_t$ and
this leads to the test $\bar\psi_\varepsilon
=1_{\left\{\left\|X_.-x_.\right\|\geq e_\alpha \right\}}$ similar to the
introduced above C-vM type test. Remind that it was shown before
that $X_t$ is the best in different senses estimator of $x_t\left(S\right)$
(see \cite{Kut94}, Section 4.3). Therefore we can think that if the
alternative is defined by the relation
$$
\left\|\frac{S\left(x_.\right)-S_0\left(x_.\right)}{\sigma
\left(x_.\right)}\right\|\geq r,
$$
then for construction of the good test statistics we have to replace
$S\left(x_t\right)$ by some estimator. Forget for instant that the Wiener
process is not differentiable and write the observed process as 
$$
\dot X_t=S\left(X_t\right)+\varepsilon \sigma \left(X_t\right)\,\dot W_t
$$
Then (formally!) $\dot X_t\rightarrow S\left(x_t\right)$ as $\varepsilon
\rightarrow 0$ and we can consider $\dot X_t$ as ``estimator'' of
$S\left(x_t\right)$. Let $\left\{\varphi _j\left(\cdot \right), j=0,\pm 1,\pm
2\ldots 
\right\}$ form an orthonormal  base in ${\cal L}_2\left(0,T\right)$, then by
Parseval identity (formally!)
we have the following equality for the corresponding (modified
$S_0\left(x_t\right)\rightarrow S_0\left(X_t\right)$) statistic 
$$
\int_{0}^{T}\left[\frac{\dot X_t-S_0\left(X_t\right)}{\sigma
\left(X_t\right)}\right]^2{\rm d}t=\varepsilon ^2\sum_{-\infty }^{\infty }
y_{j,\varepsilon }^2
$$
where the Fourier coefficients $y _{j,\varepsilon }$  we can write as follows
$$
y_{j,\varepsilon }=\int_{0}^{T}\frac{\varphi _j\left(t \right)}{\varepsilon
\sigma \left(X_t\right)}\left[\dot X_t-S_0\left(X_t\right)\right]\, {\rm d}t
=\int_{0}^{T}\frac{\varphi _j\left(t \right)}{\varepsilon \sigma
\left(X_t\right)} \left[{\rm d}X_t-S_0\left(X_t\right){\rm d}t\right].
$$
This last integral has mathematical meaning and starting from this definition
of $y_{j,\varepsilon }$ we can introduce the statistic
$$
\delta_\varepsilon
^*=\frac{1}{\sqrt{4m}}\sum_{\left|j\right|<m}^{}\left[y_{j,\varepsilon
}^2-1\right]  
$$
 Note that under ${\scr H}_0$ the random variables $y_{j,\varepsilon },
j=0,\pm 1,\pm 2\ldots$ are independent Gaussian
$$
y_{j,\varepsilon }=\int_{0}^{T}\varphi_j \left(t\right)\;{\rm d}W_t\quad \sim\quad
{\cal N}\left(0, 1\right). 
$$
Therefore the statistic $\sum_{\left|j\right|<m}^{}\left[y_{j,\varepsilon
}^2-1\right]  $ has Chi-Square distribution and the equation for $c_\alpha $
$$
\Pb_0\left\{\delta_\varepsilon ^*>c_\alpha \right\}=\alpha
$$
can be easily solved. Moreover, $m\rightarrow \infty $ we have the convergence
(under hypothesis)
$$
\frac{1}{\sqrt{4m}}\sum_{\left|j\right|<m}^{}\left[y_{j,\varepsilon
}^2-1\right] \Longrightarrow  {\cal N}\left(0,1\right).
$$

Hence we can introduce the  Ch-S test as
$$
\psi_\varepsilon^* \left(X^\varepsilon
\right)=1_{\left\{\delta_\varepsilon  ^*> z_\alpha \right\}},
$$
where $m=m\left(\varepsilon \right)\rightarrow \infty $ and $z_\alpha $ is
$1-\alpha $ quantile of the Gaussian  ${\cal N}\left(0,1 \right)$ law. 

This leads us to  the following result. 

\begin{proposition}
\label{P}
Let us suppose that $\sigma \left(x\right)^2>0, x\in R$, then the test
$\psi_\varepsilon^* \left(X^\varepsilon \right) $ belongs to ${\cal K}_\alpha
$.
\end{proposition}

We see that we need not even to use the convergence $X_t$ to $x_t$ and the
proposition is valid, say, for ergodic diffusion processes with $\varepsilon
=1$ and the asymptotic $T\rightarrow \infty $. The choice of asymptotic  is
important in the calculation of the power function. 

Suppose that the observed process has the trend coefficient $S\left(\cdot
\right)\not =S_0\left(\cdot \right)$ and the condition ${\cal R}$ is
fulfilled, then
$$
y_{j,\varepsilon }= \int_{0}^{T}\frac{\varphi _j\left(t\right)\left[S\left(X_t
\right)-S_0\left(X_t \right)\right] }{\varepsilon \sigma
\left(X_t\right)}\,{\rm d}t+\int_{0}^{T}\varphi _j\left(t\right)\,{\rm
d}W_t=z_{j,\varepsilon }+\zeta _j,
$$
where
$$
\varepsilon ^2\sum_{j}^{}z_{j,\varepsilon }^2=\int_{0}^{T}\left(\frac{S\left(X_t
\right)-S_0\left(X_t \right) }{ \sigma \left(X_t\right)}\right)^2\,{\rm
d}t\longrightarrow  \left\|\frac{S\left(x_.
\right)-S_0\left(x_. \right) }{ \sigma \left(x_.\right)} \right\|^2\geq r^2.
$$
Hence for any fixed contiguous  alternative
$S\left(x\right)=S\left(x\right)+\varepsilon 
h\left(x\right)\;\sigma \left(x\right)$ we have
\begin{align*}
\delta_\varepsilon^*&=\frac{1}{\sqrt{4m}}\sum_{\left|j\right|<m}^{}\zeta _j^2-
\frac{2}{\sqrt{4m}}\sum_{\left|j\right|<m}^{} z_{j,\varepsilon }\, \zeta
_j+\frac{1}{\sqrt{4m}}\sum_{\left|j\right|<m}^{}z_{j,\varepsilon }^2 .
\end{align*}
As $m\rightarrow \infty $ the relation 
$$
\sum_{\left|j\right|<m}^{}z_{j,\varepsilon}^2=
\left\|h\left(x_.\right)\right\|^2\left(1+o\left(1\right)\right)\geq \frac{r^2
}{\varepsilon ^2}\left(1+o\left(1\right)\right)
$$
holds. 
To have non degenerate limit 
$$
\frac{1}{\sqrt{4m}}\sum_{\left|j\right|<m}^{}z_{j,\varepsilon }^2
=\frac{r^2\,u}{\varepsilon
^2\sqrt{4m}}\left(1+o\left(1\right)\right)\longrightarrow u\geq 1
$$
we can put $m=r^4/\left(4\varepsilon ^4\right)$. Note that in this case
$$
\Ex \left(\frac{2}{\sqrt{4m}}\sum_{\left|j\right|<m}^{} z_{j,\varepsilon }\,
\zeta _j\right)^2=\frac{1}{m}\sum_{\left|j\right|<m}^{} z_{j,\varepsilon
}^2=\frac{r^2u}{\varepsilon ^2\,m}\longrightarrow 0. 
$$
Hence for the power function we obtain the limit
$$
\beta \left(\delta _\varepsilon ^*,h\right)\longrightarrow \Pb\left\{\zeta
>z_\alpha -u\right\},\qquad \quad \zeta \sim {\cal N}\left(0,1\right) .
$$ 
Another way is to fix first $m\rightarrow \infty $ such that $m\varepsilon
^{1/2}\rightarrow \infty $ and then to consider the
alternatives running to the hypothesis : $r=r\left(\varepsilon
\right)=\varepsilon \sqrt{4m}\rightarrow 0$.

Of course, this test is not uniformly consistent with the same explication as
above. 

To have uniformly consistent and minimax GoF testing we need to
control the derivatives. If we suppose that the function $h\left(\cdot
\right)$ defining the alternative is $k$ times differentiable and the ${\cal L}_2$
norm of the $k$-derivative is bounded by some constant, then we can show that
the test
$$
\delta _\varepsilon =\sum_{\left|j\right|<m}^{}w_j\,\left[y_{j,\varepsilon
}^2-1\right] 
$$ 
with weights
$$
w_i=z^2\left(1-\left|\frac{i}{m}\right|^{2k }\right),\qquad
z=\left(2\sum_{i=-m}^{m}\left[1-\left|\frac{i}{m}\right|^{2k }\right]^2 
\right)^{-1/4}
$$
and special choice of $m\rightarrow \infty $ is asymptotically minimax. The
proof is based on the approach developed by Ermakov \cite{Er} and Ingster and
Suslina \cite{IS-03} for the model {\sl signal in white Gaussian noise}. See as
well the similar problem for Poisson processes studied by  Ingster and
Kutoyants \cite{IK-07}.

\section{On local time and GoF testing}

The local time of the diffusion process \eqref{1} is defined as 
\begin{equation}
\label{ltd}
\Lambda _T\left(x\right)=\lim_{\nu \downarrow
0}\frac{\varepsilon ^2}{2\nu 
}\int_{0}^{T}1\zs{\left\{\left|X_t-x\right|\leq  \nu \right\}}\;\sigma
\left(X_t\right)^2\;{\rm d} t
\end{equation}
and admits the Tanaka-Meyer representation (see \cite{RY})
$$
\Lambda
_T\left(x\right)=\left|X_T-x\right|-\left|x_0-x\right|-\int_{0}^{T}
\sgn\left(X_t-x\right) \;{\rm d}X_t.  
$$
The local time $\Lambda _T\left(x\right)$ recently started to play an
important role in statistical inference
\cite{Kut97},\cite{BD},\cite{Kut04}. Note that in ergodic case ($\varepsilon
\equiv 1$ and $T\rightarrow \infty $), the local time is asymptotically normal
:
$$
\eta _T\left(x\right)=\sqrt{T}\left(\frac{ \Lambda _T\left(x\right)}{T\sigma
\left(x\right)^2}-f\left(x\right)\right) \Longrightarrow {\cal
N}\left(0,d_f\left(x\right)^2\right) ,
$$
where $f\left(x\right)$ is invariant density and 
$$
d_f\left(x\right)^2=4f\left(x\right)^2\;\Ex\left(\frac{1\zs{\left\{\xi
>x\right\}}-F\left(\xi \right) }{\sigma \left(\xi \right)f\left(\xi
\right)}\right)^2. 
$$
Here $\xi $ is random variable with density $f\left(x\right)$  and
$F\left(x\right)$ is its distribution function (see \cite{Kut04}, Proposition
1.25). Moreover this normed difference convergence weakly to the limit Gaussian
process in the space of continuous on ${\cal R}$ functions vanishing in
infinity \cite{Kut04}, Theorem 4.13. This property can be used in the
construction of the GoF tests as follows. Let us introduce the
 C-vM and K-S type statistics
$$
\delta _T=\int_{-\infty }^{\infty }\eta _T\left(x\right)^2{\rm
d}x ,\qquad \gamma _T=\sup_x\left|\eta _T\left(x\right)\right|
$$
 and the corresponding tests $\psi_T=1_{\left\{\delta _T>c_\alpha  \right\}}$
and $\phi_T=1_{\left\{\gamma  _T >d_\alpha \right\}}$. Then using this weak
convergence we can define the constants $c_\alpha ,d_\alpha $ (see
\cite{Kut04}, section 5.4 and Gassem \cite{Gas08}).

We can consider the similar problem in the {\it asymptotics of small noise},
i.e., to use the limit behavior of the local time in the construction of the
GoF tests.  

Let us introduce the space ${\cal L}^2\left(x_0,x_T\right)$ of square
integrable functions on $\left[x_0,x_T\right]$, where $x_t$ is solution of the
ordinary differential equation
$$
\frac{{\rm d}x_t}{{\rm d}t}=S_0\left(x_t\right),\qquad x_0,
$$
where $S_0\left(x\right)>0$. According to \eqref{ltd} we have
\begin{align*}
&\lim_{\varepsilon \rightarrow 0}\frac{\Lambda _T\left(x\right)}{\varepsilon
^2} =\lim_{\varepsilon \rightarrow 0}\lim_{\nu \rightarrow 0}\frac{1}{2\nu
}\int_{0}^{T} 1_{\left\{\left|X_t-x\right|\leq \nu \right\}}\sigma
\left(X_t\right)^2{\rm d}t\\ 
&\qquad = \lim_{\nu \rightarrow 0}\frac{1}{2\nu
}\int_{0}^{T} 1_{\left\{\left|x_t-x\right|\leq \nu \right\}}\sigma
\left(x_t\right)^2{\rm d}t=\lim_{\nu \rightarrow 0}\frac{1}{2\nu
}\int_{x_0}^{x_T} \frac{1_{\left\{\left|y-x\right|\leq \nu \right\}}\sigma
\left(y\right)^2}{S_0\left(y\right)}{\rm d}y\\ 
&\qquad =\lim_{\nu \rightarrow
0}\frac{1}{2\nu }\int_{x-\nu }^{x+\nu } \frac{\sigma
\left(y\right)^2}{S_0\left(y\right)}{\rm d}y=\frac{\sigma
\left(x\right)^2}{S_0\left(x\right)},\quad {\rm for}\quad x\in
\left[x_0,x_T\right].
\end{align*}t

We say that the random process $\eta _\varepsilon \left(x\right),x_0\leq x\leq
x_T$ converges weakly in ${\cal L}^2\left(x_0,x_T\right)$ to the random
process $\eta \left(x\right),x_0\leq x\leq x_T$ if for any function
$h\left(\cdot \right)\in {\cal L}^2\left(x_0,x_T\right)$ we have
$$
\int_{x_0}^{x_T}h\left(x\right)\,\eta _\varepsilon \left(x\right)\,{\rm
d}x\Longrightarrow  \int_{x_0}^{x_T}h\left(x\right)\,\eta  \left(x\right)\,{\rm
d}x.
$$

\begin{proposition}
\label{T3}
Let the condition ${\cal R}$ be fulfilled, then for the local time $\Lambda
_T\left(\cdot \right)  $ we have the weak convergence in ${\cal
L}^2\left(x_0,x_T\right)$ :
\begin{equation}
\label{lt}
\eta _\varepsilon \left(x\right)=\frac{1}{\varepsilon}\int_{x_0}^{x}\left( \frac{1}{S_0\left(y\right)}
-\frac{\Lambda _T\left(y\right) }{\varepsilon ^2\sigma
\left(y\right)^2}\right)\,{\rm d}y\Longrightarrow\eta  \left(x\right)=
W\left(\int_{x_0}^{x} \frac{\sigma \left(y\right)^2}{S_0\left(y\right)^3}\:{\rm
d}y\right),
\end{equation}
where $ x \in \left[x_0, x_T \right]$ and $W\left(\cdot \right)$ is a
Wiener process. 
\end{proposition}
{\bf Proof.} Note that the local time   allows us to write the equality (see
\cite{RY}) 
$$
\int_{0}^{T}h\left(X_t\right)\,{\rm d}t=\int_{-\infty }^{\infty
}h\left(x\right)\,\frac{\Lambda _T\left(x\right)}{\varepsilon ^2\sigma
\left(x\right)^2}\,{\rm d}x .
$$
Suppose that $h\left(\cdot \right)\in {\cal C}^1$. We know that
$$
\int_{0}^{T}h\left(X_t\right)\,{\rm d}t\longrightarrow
\int_{0}^{T}h\left(x_t\right)\,{\rm d}t
=\int_{0}^{T}\frac{h\left(x_t\right)}{S_0\left(x_t\right)}\,{\rm d}x_t=
\int_{x_0}^{x_T} \frac{h\left(x\right)}{S_0\left(x\right)}\,{\rm d}x
$$
and
\begin{align*}
\int_{0}^{T}\frac{h\left(X_t\right)-h\left(x_t\right)}{\varepsilon }\,{\rm
d}t&\longrightarrow \int_{0}^{T}h'\left(x_t\right)\,
x_t^{\left(1\right)}\,{\rm d}t\\
&=\int_{0}^{T}h'\left(x_t\right)\,
S_0\left(x_t\right)W\left(\int_{x_0}^{x_t}\frac{\sigma
\left(y\right)^2}{S_0\left(y\right)^3 }{\rm d}y\right) {\rm d}t\\
&=\int_{x_0}^{x_T}h'\left(x\right)\,
W\left(\int_{x_0}^{x}\frac{\sigma
\left(y\right)^2}{S_0\left(y\right)^3 }{\rm d}y\right) {\rm d}x
\end{align*}
 Hence if we denote
$$
\eta _\varepsilon \left(x\right)=\varepsilon ^{-1}\left(\frac{\Lambda
_T\left(x\right) }{\varepsilon ^2\sigma 
\left(x\right)^2}-\frac{1_{\left\{x_0\leq x\leq
x_T\right\}}}{S_0\left(x\right)}\right),\qquad x\in \RR, 
$$
and
$$
g\left(x\right)=\int_{x_0}^{x}
\frac{\sigma\left(y\right)^2}{S_0\left(y\right)^3 }{\rm d}y,\qquad x\in
\left[x_0,x_T\right] 
$$
then we can write
\begin{equation*}
\label{lt-1}
\int_{-\infty }^{\infty }h\left(x\right)\,\eta _\varepsilon \left(x\right)\,{\rm
d}x\longrightarrow  
\int_{x_0}^{x_T}h'\left(x\right)\,W\left(g\left(x\right) \right)\;{\rm d}x.
\end{equation*}
The same time integrating by parts we have
$$
\int_{-\infty }^{\infty }h\left(x\right)\,\eta _\varepsilon
\left(x\right)\,{\rm
d}x=\int_{-\infty}^{\infty}h'\left(x\right)\,\psi_\varepsilon \left(x
\right)\;{\rm d}x,
$$
where we put
$$
\psi_\varepsilon \left(x \right)=\int_{x}^{x_0}\eta _\varepsilon \left(y\right)\,{\rm
d}y.
$$

Hence
\begin{equation*}
\label{lt-2}
\int_{-\infty }^{\infty }h'\left(x\right)\,\psi_\varepsilon \left(x \right)\;{\rm
d}x\longrightarrow  \int_{x_0}^{x_T}h'\left(x\right)\,W \left(g\left(x
\right)\right)\;{\rm d}x.
\end{equation*}
We see that the values of $h\left(\cdot \right)$ outside of the interval
$\left[x_0,x_T\right]$ have no contribution in the limit. Therefore we have as
well the convergence
\begin{equation*}
\label{lt-3}
\int_{x_0}^{x_T }h'\left(x\right)\,\psi_\varepsilon \left(x \right)\;{\rm
d}x\longrightarrow  \int_{x_0}^{x_T}h'\left(x\right)\,W \left(g\left(x
\right)\right)\;{\rm d}x.
\end{equation*}
Remind that this is true for any function $h\left(\cdot \right)\in {\cal
C}^1$. Therefore the proposition is proved.

\bigskip

The convergence \eqref{lt} suggests the construction of the following
test. Let 
$$
\delta _\varepsilon =\left(\varepsilon \int_{x_0}^{x_T}\frac{\sigma
\left(x\right)^2}{S_0\left(x\right)^3}{\rm
d}x\right)^{-2}\int_{x_0}^{x_T}\frac{\sigma
\left(x\right)^2}{S_0\left(x\right)^3}\left(\int_{x_0}^{x}
\left(\frac{1}{S_0\left(y\right)}- \frac{\Lambda _T\left(y\right)}{\varepsilon
^2\sigma \left(x\right)^2}\right){\rm d}y\right)^2  {\rm d}x . 
$$
Then it can be shown that 
$$
\delta _\varepsilon \Longrightarrow \int_{0}^{1}w_v^2\,{\rm d}v
$$
and the corresponding test $\psi_\varepsilon =1_{\left\{\delta _\varepsilon
>c_\alpha \right\}}$ is asymptotically distribution free. 

\bigskip

{\bf Remark.} We see that despite the ergodic case, the local time random
function has no limit (as process) and for small values of $\varepsilon $ its
behavior is close to white noise process.

\bigskip

{\bf Remark.} We supposed above that $S_0\left(x\right)>0$ for all $x\in
\left[x_0,x_T\right]$. In the case $S_0\left(x_0\right)=0$ the deterministic
solution  $x_t\equiv x_0$ and we have the following basic hypothesis  
$$
{\scr H}_0\quad :\qquad {\rm d}X_t=\varepsilon \sigma \left(X_t\right)\;{\rm
d}W_t,\qquad X_0,\quad 0\leq t\leq T. 
$$	
The test can be based on the statistic
$$
\delta _\varepsilon =\int_{0}^{T}\left(\frac{X_t-x_0}{T\varepsilon\,\sigma
\left(x_0\right) }\right)^2{\rm d}t 
$$
and it is easy to see that (under hypothesis ${\scr H}_0\ $)
$$
\delta _\varepsilon \longrightarrow \int_{0}^{1}w_s^2\;{\rm d}s.
$$
If for some $x_*>x_0$ we have  $S_0\left(x\right)>0, x\in[x_0,x_*)$ and
$S_0\left(x_*\right)=0$, then by Lipschitz condition
$$
t=\int_{x_0}^{x_t}\frac{{\rm d}x}{S_0\left(x\right)-S_0\left(x_*\right) }\geq 
\frac{1}{L}\int_{x_0}^{x_t}\frac{{\rm
d}x}{x_*-x}=\frac{1}{L}\ln\frac{x_*-x_0}{x_*-x_t} 
$$
and we see that the equality $x_t=x_*$ is impossible (well known property).

\section{On composite basic hypothesis}

Suppose that under hypothesis ${\scr H}_0$ the observed diffusion process is
solution of the stochastic differential equation
$$
{\rm d}X_t=S\left(\vartheta ,X_t\right)\,{\rm d}t+\varepsilon \sigma
\left(X_t\right)\,{\rm d}W_t,\quad X_0=x_0,\quad 0\leq t\leq T, 
$$
where the trend coefficient $S\left(\vartheta ,x\right)$ is a known function
which depends on unknown parameter $\vartheta \in \Theta
=\left(\beta,\gamma  \right)$. Then the limit solution $x_t$ depends on the
true value $\vartheta $, i.e., $x_t=x_t\left(\vartheta
\right)$ and the natural modification of the test statistic can be based on
the normalized difference 
$$
Y_{t }(\hat\vartheta_\varepsilon)=\frac{X_t-x_t(\hat\vartheta
_\varepsilon )}{\varepsilon } ,
$$
where as $\hat\vartheta
_\varepsilon$ we can take, say, the maximum likelihood estimator. We suppose
that the functions $S\left({\vartheta ,x}\right)$ and $\sigma \left(x\right)$ are positive
and sufficiently smooth to calculate the derivatives below and to provide the
``usual properties of estimators''.  Remind that
(under regularity conditions) this estimator is consistent and asymptotically
normal (see \cite{Kut94}, Theorem 2.2). Moreover the MLE admits the
representation (see \cite{Kut94}, Theorem 3.1)
$$
\hat\vartheta_\varepsilon=\vartheta _0+\varepsilon {\rm I}_T\left(\vartheta
_0\right)^{-1} \int_{0}^{T}\frac{\dot S\left(\vartheta _0,x_t\left(\vartheta
_0\right)\right)}{\sigma \left(x_t\left(\vartheta _0\right)\right)}\;{\rm
d}W_t+o\left(\varepsilon \right) ,
$$
where $ \vartheta _0$ is the true value of $\vartheta $ and ${\rm
I}_T\left(\vartheta_0 \right)$ is the Fisher information:
$$
{\rm I}_T\left(\vartheta \right)=\int_{0}^{T}\left(\frac{\dot S\left(\vartheta
,x_t\left(\vartheta\right)\right)}{\sigma (x_t\left(\vartheta\right))}
\right)^2{\rm d}t.
$$
 It can be shown (see Rabhi \cite{Rab}), that in regular (smooth) case the
limit distribution of 
$$
\delta _\varepsilon=\int_{0}^{T}Y_{t
}(\hat\vartheta_\varepsilon)^2\,{\rm d}t
 $$ 
coincides with the distribution of the following integral
$$ 
\int_{0}^{T}\left(x_t^{\left(1\right)}\left(\vartheta _0\right)-{\rm
I}_T\left(\vartheta_0 \right)^{-1}\dot x_t\left(\vartheta _0\right)
\int_{0}^{T}\frac{\dot S\left(\vartheta _0,x_s\left(\vartheta
_0\right)\right)}{\sigma \left(x_s\left(\vartheta _0\right)\right)} {\rm d}W_s
\right)^2{\rm d}t,
$$
 where $ \dot x_t\left(\vartheta \right)=\frac{\partial }{\partial \vartheta
}\;x_t\left(\vartheta \right) $ and $x_t^{\left(1\right)}\left(\vartheta
\right)$ is solution of the linear equation
$$
{\rm d}x_t^{\left(1\right)}=S'\left(\vartheta
,x_t\right)x_t^{\left(1\right)}\,{\rm d}t+ \sigma
\left(x_t\right)\,{\rm d}W_t,\quad x_0^{\left(1\right)}= 0,\quad 0\leq t\leq T, 
$$
The test is no more distribution free, but in some cases it can be done
asymptotically distribution free  (ADF) if the  second limit $T\rightarrow \infty $
is taken \cite{Rab}. 

Another possibility to have an ADF test is to use {\it the estimator process }
$\hat \vartheta _{t,\varepsilon }, 0\leq t\leq T$, where $\hat \vartheta
_{t,\varepsilon }$ is an estimator   constructed by the observations
$X^t=\left\{X_s,0\leq s\leq t\right\}$. Of course, we have to suppose that
this estimator is consistent and asymptotically normal for all values of $t\in
(0,T]$. For example, in the   linear case
$$
{\rm d}x_t=\vartheta \,h\left(X_t\right){\rm d}t+\varepsilon \sigma
\left(X_t\right)\,{\rm d}W_t, \quad X_0=x_0, \quad 0\leq t\leq T
$$
we can take the MLE process
$$
\hat\vartheta _{t,\varepsilon
}=\left(\int_{0}^{t}\frac{h\left(X_s\right)^2}{\sigma
\left(X_s\right)^2}\;{\rm d}s\right)^{-1}
\int_{0}^{t}\frac{h\left(X_s\right)}{\sigma \left(X_s\right)^2}\;{\rm d}X_s,
\quad 0< t\leq T.
$$
Sometimes an estimator process can have recurrent structure (see Levanoy {\it
  et al.} \cite{LSZ}).  If $\hat \vartheta _{t,\varepsilon }$ is the MLE, then
  it can be shown (see the similar calculus above \eqref{cal})
  the limit of $\delta _\varepsilon $ is (below $x_s=x_s\left(\vartheta _0\right)$)
\begin{align*}
&\int_{0}^{T}S\left(\vartheta _0,x_t\right)^2 \left(\int_{0}^{t}\frac{\sigma
\left(x_s\right)}{S\left(\vartheta _0, x_s\right)}{\rm d}W_s\right.\\ &\qquad
\qquad \left. -{\rm I}_t\left(\vartheta_0 \right)^{-1}\int_{0}^{t}\frac{\dot
S\left(\vartheta _0,x_s\right)}{S\left(\vartheta _0,x_s\right)}{\rm d}s
\int_{0}^{t}\frac{\dot S\left(\vartheta _0,x_s\right)}{\sigma
\left(x_s\right)} {\rm d}W_s \right)^2{\rm d}t\\ &\qquad \equiv
\int_{0}^{T}S\left(\vartheta _0,x_t\right)^2Z\left(t,\vartheta _0\right)^2{\rm
d}t.
\end{align*}
The process $Z\left(t,\vartheta _0\right)$ has stochastic differential
$$
{\rm d}Z\left(t,\vartheta _0\right)=A \left(t,\vartheta _0\right){\rm
d}t+B\left(t,\vartheta _0\right){\rm d}W_t,\qquad Z\left(0,\vartheta
_0\right)=0. 
$$
with the corresponding random function $A \left(t,\vartheta _0\right)$ and
deterministic $B\left(t,\vartheta _0\right)$. Hence
$$
B\left(t,\vartheta _0\right)^{-1} \left(Z\left(t,\vartheta
_0\right)-\int_{0}^{t}A \left(s,\vartheta _0\right){\rm d}s\right)=W_t 
$$
and
\begin{align*}
\delta _0 =\int_{0}^{T} \left(\frac{Z\left(t,\vartheta 
_0\right)-\int_{0}^{t}A \left(s,\vartheta _0\right){\rm d}s}{T\,S\left(\vartheta
_0,x_t\right)B\left(t,\vartheta _0\right) }\right)^2{\rm d}t=\int_{0}^{1}w_v^2\;{\rm d}v
\end{align*}

Using ``empirical versions'' of these functions  it is possible to construct
 an ADF test based on the following statistics 
$$
\delta _\varepsilon  =\int_{\varepsilon^\mu  }^{T} \left(\frac{Y_t
(\hat\vartheta  
_{t,\varepsilon} )-\int_{0}^{t}A_\varepsilon  (s,\hat\vartheta 
_\varepsilon ){\rm d}s}{T\,S(\hat\vartheta 
_\varepsilon ,X_t)B_\varepsilon (t,\hat\vartheta 
_\varepsilon ) }\right)^2{\rm d}t\Longrightarrow \int_{0}^{1}w_v^2\;{\rm d}v.
$$
What is empirical version of stochastic integral we explain below. The
integral is started at $t=\varepsilon ^\mu $, where $\mu \in \left(0,1\right)$
because for the values $t\in \left[0,\varepsilon \right]$ the  estimator
$\hat\vartheta _{t,\varepsilon }$ is not asymptotically normal. 

 One else ADF  test can be constructed by ``compensating'' the
 additional random part by the following way. First we rewrite the stochastic
 integral  (It\^o formula)
\begin{align*}
H_\varepsilon \left(\vartheta \right)& =\int_{0}^{T}\frac{\dot S\left(\vartheta
,X_s\right)}{\sigma \left(X_s\right)^2}  \left[{\rm d}X_s-S\left(\vartheta
,X_s\right){\rm d}s \right]   =    \int_{0}^{T}\frac{\dot
S\left(\vartheta,X_s\right)}{\sigma \left(X_s\right)}  {\rm d}W_s 
\\
&=\int_{x_0}^{X_T}\frac{\dot S\left(\vartheta,y\right) }{\sigma
\left(y\right)^2}  {\rm d}y- \varepsilon ^2\int_{0}^{T}\frac{\dot
S'\left(\vartheta,X_s\right)\sigma \left(X_s\right) -2\dot
S\left(\vartheta,X_s\right)\sigma' \left(X_x\right)}{2\sigma \left(X_s\right)
}{\rm d}t\\
&\quad - \int_{0}^{T}\frac{\dot S\left(\vartheta
,X_s\right)\,S\left(\vartheta ,X_s\right)}{\sigma \left(X_s\right)^2}  {\rm d}s.
\end{align*}
The last expression for $H_\varepsilon \left(\vartheta \right)$ does not
contain stochastic integral and we can put the estimator, i.e., the random
variable $H_\varepsilon \left(\hat\vartheta_\varepsilon  \right)$ is well
defined (empirical version).
 Then introduce the stochastic process
\begin{align*}
Y _\varepsilon(t,\hat\vartheta _\varepsilon ) =\frac{X_t-x_t(\hat\vartheta
_\varepsilon )}{\varepsilon } + {\rm I}_T(\hat\vartheta _\varepsilon )^{-1} \dot
x_t(\hat\vartheta _\varepsilon ) H_\varepsilon (\hat\vartheta _\varepsilon ) .
\end{align*}
Note that  it can be easily shown that
\begin{align*}
\sup_{0\leq t\leq T}\left| \dot x_t(\hat\vartheta _\varepsilon )-\dot
x_t\left(\vartheta _0 \right)\right|\rightarrow 0 ,\quad {\rm I}_T(\hat\vartheta
_\varepsilon )\rightarrow {\rm I}_T\left(\vartheta _0 \right),\quad
H_\varepsilon (\hat\vartheta _\varepsilon )\rightarrow H_0 \left(\vartheta
_0\right),
\end{align*}
where
$$
H_0 \left(\vartheta_0\right)=\int_{0}^{T}\frac{\dot
S\left(\vartheta _0,x_s\left(\vartheta _0\right)\right)}{\sigma
\left(x_s\left(\vartheta _0\right)\right)} {\rm d}W_s .
$$
Hence the stochastic process $Y_\varepsilon (t,\hat\vartheta _\varepsilon
)$ converges uniformly on $t\in\left[0,T\right]$ to the Gaussian process
$x_t^{\left(1\right)}\left(\vartheta _0\right)$ and we can use the statistic
$$
\delta _\varepsilon =\left[\int_{0}^{T}\frac{\sigma
\left(X_t\right)^{2}}{S(\hat\vartheta_\varepsilon ,X_t)^{2}}{\rm
d}t\right]^{-2} \; \int_{0}^{T} \frac{Y_\varepsilon (t,\hat\vartheta
_\varepsilon ) ^2\,\sigma \left(X_t\right)^2}{S(\hat\vartheta
_\varepsilon ,X_t)^4} {\rm d}t\Longrightarrow \int_{0}^{1}w_v^2\,{\rm d}v.
$$
We see that the test $\hat\psi _\varepsilon =1_{\left\{\delta
_\varepsilon>c_\alpha  \right\}} $ based on this statistics is ADF.

This test has to be consistent against any fixed alternative. Indeed, let the
observed process be
$$
{\rm d}X_t=S\left(X_t\right)\,{\rm d}t+\varepsilon \sigma
\left(X_t\right)\;{\rm d}W_t , \quad X_0=x_0,\quad 0\leq t\leq T
$$
and the limit solution $x_t\left(S\right)$ satisfies 
$$
g=\inf_{\vartheta \in \Theta }\left\|x_.\left(S\right)-x_.\left(\vartheta
\right)\right\| >0.
$$
The MLE $\hat\vartheta _\varepsilon $ in this misspecified situation converges
to the value $\vartheta _*$ which minimizes the Kullback-Leibner distance (see
\cite{Kut94}, Section 2.6) 
$$
\vartheta _*=\arg\inf _{\vartheta \in \Theta }
\int_{0}^{T}\left(\frac{S\left(\vartheta
,x_t\left(S\right)\right)-S\left(x_t\left(S\right)\right)}{\sigma
\left(x_t\left(S\right)\right)} \right)^2{\rm d}t.
$$
Hence
\begin{align*}
Y_\varepsilon \left(t,\hat\vartheta _\varepsilon
\right)=x_t^{\left(1\right)}\left(S\right) +
\frac{x_t\left(S\right)- x_t\left(\vartheta _*\right)}{\varepsilon }+{\rm
I}_T\left(\vartheta _*\right)^{-1} \dot x_t\left(\vartheta
_*\right)H_0\left(\vartheta _*\right) +o\left(1\right)
\end{align*}
and
$$
\delta _\varepsilon \longrightarrow \infty . 
$$
Therefore the test is consistent.

\end{document}